\documentclass[conference]{IEEEtran}
\IEEEoverridecommandlockouts

\usepackage{graphicx} 
\usepackage{amsmath}
\usepackage{amsthm}
\usepackage{amssymb}
\usepackage{xcolor}
\usepackage{verbatim}
\usepackage{algorithm}
\usepackage[noend]{algpseudocode}
\usepackage[bottom=2.0 cm,top=2.0 cm, left=2.0 cm, right=2.0 cm]{geometry}
\usepackage{bbm}
\usepackage[bb=boondox]{mathalfa}
\usepackage{graphicx} %
\usepackage{hyperref}
\hypersetup{colorlinks=true,
            linkcolor=black,
            anchorcolor=black,
            citecolor=black}

\newcommand{\A}{\mathbb{A}}

\newcommand{\B}{\mathbb{B}}

\newcommand{\C}{\mathbb{C}}

\newcommand{\Sub}{\mathbf{S}}
\newcommand{\T}{\mathbf{T}}

\newcommand{\Z}{\mathbf{Z}}
\newcommand{\X}{\mathbf{X}}
\newcommand{\J}{\mathbb{J}}

\newcommand{\Y}{\mathbf{Y}}
\newcommand{\D}{\mathbf{D}}
\newcommand{\I}{\mathbb{I}}
\newcommand{\0}{\mathbb{0}}
\newcommand{\PL}{\text{\L}}
\newcommand{\prox}{\texttt{prox}}
\definecolor{xkchen}{rgb}{0.75, 0.2, 1.0}
\newtheorem{assumption}{Assumption}

\newtheorem{theorem}{Theorem}

\newtheorem{lemma}{Lemma}

\usepackage{cite}
\usepackage{amsmath,amssymb,amsfonts}
\usepackage{graphicx}
\usepackage{textcomp}
\usepackage{xcolor}
\def\BibTeX{{\rm B\kern-.05em{\sc i\kern-.025em b}\kern-.08em
    T\kern-.1667em\lower.7ex\hbox{E}\kern-.125emX}}
\begin{document}

\title{A Parameter-free Decentralized Algorithm for Composite Convex Optimization
\thanks{Chen and Scutari are with the School of Industrial Engineering, Purdue University; emails:\texttt{<chen4373,gscutari>@purdue.edu}.  Kuruzov and Gasnikov are with  the Innopolis University; emails:  \texttt{kuruzov.ia@phystech.edu, gasnikov@yandex.ru}. }
}
\author{Xiaokai Chen, Ilya Kuruzov, Gesualdo Scutari, and Alexander Gasnikov }\vspace{-1.5cm}

\maketitle
\begin{abstract}
The paper studies decentralized optimization over networks, where agents minimize a composite objective consisting of the sum of smooth convex functions--the agents' losses--and an additional nonsmooth convex extended value function.   
We propose a  decentralized algorithm wherein  agents {\it adaptively} adjust their stepsize using local backtracking procedures that require {\it no global} (network) information or extensive inter-agent communications. 
Our adaptive decentralized method enjoys robust convergence guarantees, outperforming  existing decentralized methods, which are not adaptive. Our design is centered on a three-operator splitting, applied to a reformulation of the  optimization problem. This reformulation utilizes a proposed  BCV metric, which facilitates decentralized implementation and local stepsize adjustments while  guarantying   convergence.\vspace{-0.2cm} 
\end{abstract}


\section{Introduction}
 Consider the following decentralized optimization  
\begin{equation}
    \label{eq:problem}
    \tag{P}
    \begin{aligned}    
&\min_{x\in\mathbb{R}^d}&\quad&u(x):=f(x)+r(x),    \end{aligned}
\end{equation}
\begin{equation}
f(x):=\frac{1}{m}\sum_{i=1}^mf_i(x),\vspace{-0.2cm}
\end{equation}
where  $f_{i} : \mathbb{R}^{d} \to \mathbb{R}$ \ is the cost function of agent $i\in [m]:=\{1,\ldots, m\}$, assumed to be convex, \emph{locally} smooth,   and known only to the agent;\ $r : \mathbb{R}^{d} \to \mathbb{R}\cup \{-\infty,+\infty\}$ \ is a convex, nonsmooth  function, known to all agents, which can be used to enforce shared constraints or specific structures on the solution (e.g., sparsity or low-rank). Agents are embedded in a {\it mesh} communication network, modeled as a fixed, undirected and connected graph $\mathcal{G}$ with no servers.

 Problem~\eqref{eq:problem} arises from several applications of interest, including signal processing, machine learning, multiagent systems, and communications. The literature abounds with decentralized solution methods for~\eqref{eq:problem} (when $f_i$'s are {\it globally} smooth); we refer to the recent tutorial (and reference therein) \cite{Nedic_Olshevsky_Rabbat2018} and book \cite{Sayed-book} 
 for a detailed review. These methods commonly rely on conservative stepsize bounds for convergence, dependent on parameters such as the Lipschitz constants of the agents’ gradients, the spectral gap of the graph adjacency matrix, or other topological properties. However, this crucial information is not locally available to the agents. Consequently, stepsize values are often set through manual tuning, leading to performance that is unpredictable, problem-dependent, and not reproducible. Furthermore, these methods fails when agents'  losses are only \emph{locally} smooth--refer to Sec.\ref{sec:related-work} for more information.  This paper brings adaptivity to the decentralized setting~\eqref{eq:problem}.

\subsection{Main contributions}

\textit{1. Algorithm design:} We introduce an adaptive decentralized algorithm tailored to the composite structure of Problem~\eqref{eq:problem}. To deal with the nonsmooth term, our design  leverages a three-operator splitting decomposition. This is applied to a reformulation of~\eqref{eq:problem} based on  the  BCV technique~\cite{bertsekas2016nonlinear,o2020equivalence} (named after Bertsekas, O'Connor, and Vandenberghe). A suitable  metric is chosen  in the BCV transformation, to facilitate  a fully decentralized implementation of the splitting-based algorithm that uses only {\it neighboring communications} while enabling  adaptive  updates of the agents' stepsize through local backtracking procedures. Remarkably, our algorithm does not require agents to know any global optimization or network parameters.  

\textit{2. Convergence guarantees:} We rigorously prove convergence of the proposed method, demonstrating a sublinear $\mathcal{O}(1/k)$ rate for a suitably defined optimality gap, where $k$ denotes the iteration count.  Preliminary numerical results highlight that our adaptive algorithm significantly outperforms existing decentralized algorithms applicable to~\eqref{eq:problem}, all of which rely on nonadaptive stepsize choices. 

\subsection{Related works}\label{sec:related-work}
\textit{1. Adaptive centralized methods:} In recent years, there has been an increasing interest in developing adaptive (a.k.a. parameter-free) optimization algorithms for centralized and federated settings. These methods include established techniques such as line-search approaches  \cite{bertsekas2016nonlinear}, Polyak's stepsize \cite{polyak1969minimization}, Barzilai-Borwein's stepsize \cite{BarzilaiBorwein1988}, as well as more recent advancements that   estimate the local curvature of the cost function \cite{Malitsky2019AdaptiveGD,Malitsky_Mishchenko_2024,Zhou24}.   Additionally, other adaptive gradient methods tailored specifically to machine learning challenges include   AdaGrad \cite{duchi2011adaptive}, Adam \cite{kingma2014adam}, AMSGrad \cite{Reddi2018OnTC}, NSGD-M \cite{pmlr-v119-cutkosky20b}, and their variants \cite{Orabona19,Ward20}. 
Some of these  algorithms have been adapted for federated architectures (server-client systems) \cite{50448,Li2022OnDA,Chen_Li_Li_2020}. However, these methods are not suitable for mesh networks because they rely  on a central server to aggregate local model updates and tune the common stepsize--a process impractical in decentralized mesh setups. 
Among these methods, only~\cite{Zhou24,polyak1969minimization,Malitsky_Mishchenko_2024,duchi2011adaptive} can deal with nonsmooth convex terms, as in Problem~\eqref{eq:problem}. 

{\it 2. Adaptive decentralized   methods:} The landscape of parameter-free {\it decentralized} methods is limited, with   notable proposals including \cite{Nazari_Tarzanagh_Michailidis_2022,chen2023convergence, li2024problem,kuruzov2024achieving,Zhou24,Notarstefano24}. Specifically, \cite{Nazari_Tarzanagh_Michailidis_2022,chen2023convergence,li2024problem} study stepsize adaptivity for stochastic (non)convex/ online and {\it smooth} optimization, typically using gradient normalization from past iterates.  Except for \cite{li2024problem}, these methods assume globally Lipschitz continuous losses, enabling convergence with standard stepsizes of the order  $\mathcal{O}(1/\sqrt{k})$ (here, $k$ is the iteration index). Additionally,   \cite{Nazari_Tarzanagh_Michailidis_2022, chen2023convergence} require prior knowledge of certain problem-dependent parameters. The work~\cite{Notarstefano24} introduces a Port-Hamiltonian System framework  for smooth, strongly convex, and unconstrained problems. While parameter-free convergence holds in the centralized setting,   convergence (global asymptotic stability) in decentralized settings is ensured  only under specific graph structures or graph-dependent stepsize constraints, which rely on global network information.       No explicit convergence rate expression is provided.  Recently, \cite{kuruzov2024achieving} proposed a parameter-free decentralized method for smooth, strongly convex instances of~\eqref{eq:problem}, establishing linear convergence and demonstrating superior theoretical and practical performance compared to  non-adaptive  decentralized algorithms.
  

None of the decentralized methods described above   handle \textit{composite} optimization problems as   \eqref{eq:problem}, nor do they account for \textit{locally} smooth losses. While proximal adaptive methods for composite optimization have been studied in centralized settings \cite{Yura-Pock-LSPrimal-dual-18,Latafat_23b}, their  adaptations to decentralized problems exhibit significant practical limitations.
   Specifically,    the decentralization of   \cite{Yura-Pock-LSPrimal-dual-18} requires at  {\it every }  line-search step    {\bf (i)}    transmission of {\it ambient-sized} vectors--causing   prohibitive communication costs--and {\bf (ii)} computation of the smallest agents' line-search  stepsize across the {\it entire} network.  The method from \cite{Latafat_23b}, while avoiding line-search procedures, still necessitates computing and propagating \textit{global} scalars at each iteration for stepsize and relies on knowledge of a {\it global network quantity} unavailable locally. 
  Resorting to conservative local overestimates of this parameter reduces adaptivity, slows convergence, and makes the algorithm sensitive to critical tuning parameters. These practical drawbacks are confirmed by numerical results in Sec.~\ref{sec:simulations}. 

 In summary, existing adaptive decentralized methods cannot solve composite optimization problems~\eqref{eq:problem} {\it using only neighbor communications}. To the best of our knowledge, the parameter-free decentralized algorithm presented here is the first to address this critical challenge.

\section{Algorithm Design}  
We study~\eqref{eq:problem} under the following   assumptions.
\begin{assumption} 
\label{ass:function}
(\textbf{i}) Each $f_i:\mathbb{R}^d\rightarrow \mathbb{R}$ is convex on $\mathbb{R}^d$ and   {{\it locally}} smooth; 
(\textbf{ii}) $r:\mathbb{R}^d\rightarrow\mathbb{R}\cup\{-\infty,\infty\}$ is convex, proper and lower-semicontinuous; and 
(\textbf{iii}) $u:\mathbb{R}^d\rightarrow \mathbb{R}$ is lower bounded. 
\end{assumption}
\begin{assumption}
\label{ass:graph}The network is modeled as an undirected, connected graph $\mathcal G=(\mathcal{V},\mathcal{E})$, where $\mathcal{V}=[m]$ and $(i,j)\in\mathcal{E}$ if and only if there is a   link between agent $i$ and $j$.  
\end{assumption}

\begin{assumption}
    [Gossip  matrices]
    \label{ass:W}
    Let $W_{\mathcal{G}}$ denote the set of symmetric, doubly stochastic, gossip matrices $\widetilde{W}:=(\widetilde{w}_{ij})_{i,j=1}^m$ that are compliant with   $\mathcal{G}$, i.e.,   $\widetilde{w}_{ii}>0$, for all $i\in[m]$;   $\widetilde{w}_{ij}>0$ for all $(i,j)\in\mathcal{E}$; and $\widetilde{w}_{ij}=0$ otherwise. 
\end{assumption}
  Assumption~\ref{ass:function}   requires only {\it local} smoothness of the losses, which significantly enlarges the class of covered functions and applications. Examples include  
covariance estimation problems and   Poisson inverse problems~\cite{bertero2009image,boyd2011distributed,bauschke2017descent} arising  from image processing.  Assumption~\ref{ass:W} is standard in the literature of gossip-based   algorithms--see, e.g., \cite{Nedic_Olshevsky_Rabbat2018,Sayed-book}.

\subsection{A three-operator splitting-based algorithm} The starting point of  our new algorithm design is the following reformulation of \eqref{eq:problem} {\it leveraging the BCV technique}: introducing local copies $x_i\in \mathbb{R}^d$ of  $x$ and slack variables $\tilde x_i\in \mathbb{R}^d$, along with the row-wise stack matrices $\mathbf{X}:=[x_1, \ldots, x_m]^\top$ and $\tilde{\mathbf{X}}:=[\tilde{x}_1, \ldots, \tilde {x}_m]^\top$, the following optimization problem
\begin{equation}
    \label{eq:ausiliary}
    \tag{P$^{\prime}$}
    \min_{\X,\widetilde{\X}\in\mathbb{R}^{m\times d}} \underbrace{F(\X)}_{:=\widetilde{F}(\X,\widetilde{\X})}+\underbrace{R(\X)+\delta_{\{\mathbf{0}\}}(\widetilde{\X})}_{:=\widetilde{R}(\X,\widetilde{\X})}+\underbrace{\delta_{\{\mathbf{0}\}}(\PL \X+M\widetilde{\X})}_{:=\widetilde{G}(\X,\widetilde{\X})}.
\end{equation}
is equivalent to the  original Problem~\eqref{eq:problem}. In \eqref{eq:ausiliary},   we defined 
 $$F(\X):=\sum_{i=1}^m f_i(x_i),\quad  R(\X):=\sum_{i=1}^m r(x_i);$$  $\delta_{\{\mathbf{0}\}}:\mathbb{R}^{m\times d}\to \mathbb{R}\cup \{\infty\}$ is the indicator function of  $\{\mathbf{0}\}$;  $\PL\in\mathbb{S}^{m}$   satisfies  $\texttt{null}(\PL)=\texttt{span}(\mathbf{1}_m)$,  and $M\in\mathbb{S}^{m}_{++}$ is the BCV metric.  The matrix  $\PL$    enforces consensus among the agents' variables $x_i$'s via  $\PL \mathbf X=\mathbf 0$.  The choice of $M$    will be shown to be  crucial to enable   a fully decentralized implementation of the proposed algorithm as well as  adaptive local stepsize selection.  This justifies the presence of $\widetilde G$ in (\ref{eq:ausiliary}), 
 which is unconvential  in the classical decentralized optimization literature.  Here, $\mathbb{S}^{m}$ (resp. $\mathbb{S}^{m}_{++}$)   denotes the set of $m\times m$ real symmetric (resp. positive definite)  matrices.

Problem~(\ref{eq:ausiliary}) is equivalent to solving the   
monotone inclusion based on the following three-operator splitting~\cite{ryu2022large}:  
\begin{equation}
    \label{eq:inclusion 2}
  \text{find }  \mathbf{T}:=\begin{bmatrix}
        \X\\\widetilde{\X}
    \end{bmatrix}\in \mathbb{R}^{2m\times d} \text{ s.t. }  0\in({\A}+{\B}+{\C})  \mathbf{T},
\end{equation}
where
$ {\A}:=\partial\widetilde{G}$,     ${\B}:=\partial \widetilde{R}$ and  $ {\C}:=\nabla \widetilde F$. 
 Invoking the  Davis-Yin three-operator splitting~\cite{davis2017three},  (\ref{eq:inclusion 2}) is equivalent to  
 \begin{equation}  \label{eq:inclusion 3}\text{find } \mathbf T, \mathbf Z\in \mathbb{R}^{2m\times d}  \,\text{ s.t. } \,\mathbb{D}_{\alpha} \mathbf Z=\mathbf Z \,\text{ and } \,\mathbf T=   \mathbb{J}_{\alpha\B} \mathbf Z,\end{equation}
where\vspace{-0.2cm}\begin{equation}\label{eq:D_op}\mathbb{D}_{\alpha}:=\I-\mathbb{J}_{\alpha \B}+\mathbb{J}_{\alpha \A}\circ(2\J_{\alpha\B}-\I-\alpha\C\circ\mathbb{J}_{\alpha\B}),\end{equation}   $\J_\A:=(\I+\A)^{-1}$ 
 is the resolvent of a given operator $\A$; and the constant $\alpha>0$  plays the role of the stepsize.  Notice that  $\J_{\alpha \partial \widetilde R}$ (resp. $\J_{\alpha \partial \widetilde G}$)   is equivalent to the  proximal operator $\texttt{prox}_{\alpha\widetilde{R}}$ (resp. $\texttt{prox}_{\alpha\widetilde{G}}$ ) applied row-wise to its argument.  

 One can solve    (\ref{eq:inclusion 3}) through the standard fixed-point  Krasnosel'skii-Mann  iteration \cite{davis2017three}: given $\Z^k$ and $\alpha^k>0$,  
\begin{equation}
\label{eq:FPI_1}
      \begin{aligned}
    &\mathbf{T}_{\B}^{k}= \texttt{prox}_{\alpha^k\widetilde{R}}\Z^k,\\
&\mathbf{T}_{\A}^{k+1}=\texttt{prox}_{\alpha^k\widetilde{G}}\left(2\mathbf{T}_{\B}^{k}-\Z^k -\alpha^k \nabla \widetilde F (\mathbf{T}_{\B}^{k})\right),\\
& \Z^{k+1}=\Z^{k}+  \mathbf{T}_{\A}^{k+1}-\mathbf{T}_{\B}^{k}.
\end{aligned} 
 \end{equation}

Next, we   develop~\eqref{eq:FPI_1} further to include a backtracking procedure that allows for adaptive selection of the stepsize.

 This development necessitates defining a suitable search direction and a proper merit function.  The backtracking procedure and the direction are established based on the saddle-point reformulation of Problem~\eqref{eq:ausiliary}: $\min_{\T\in\mathbb{R}^{2m\times d}}\max_{\Sub\in\mathbb{R}^{2m\times d}}\mathcal{L}(\T,\Sub)$, where  
\begin{equation}
    \label{eq:new Lagarangian}
    \mathcal{L}(\T,\Sub):=\widetilde{F}(\T)+\widetilde{G}(\T)+\langle\Sub,\T\rangle-\widetilde{R}^*(\Sub),
\end{equation}
  is the Lagrangian function,   with $\T$ as the primal variables and $\Sub$ as the dual variables, and $\widetilde{R}^*$ denoting the conjugate of $\widetilde{R}$.
Introducing the intermediate variable 
$$\Sub^k:=\frac{1}{\alpha^k}(\Z^k-\mathbf{\T}_{\B}^k),$$ allows us to express $\T_{\A}^{k+1}$ in~\eqref{eq:FPI_1} as  \begin{equation}
    \label{eq:primal update}
\T_{\A}^{k+1}=\texttt{prox}_{\alpha^k,\widetilde{G}}(\T_{\B}^k-\alpha^k\Sub^k-\alpha^k\nabla \widetilde{F}(\T_{\B}^k)).
\end{equation}
If the stepsize $\alpha^k>0$ satisfies the following backtracking condition:  given  $\delta\in(0,1)$, find the largest $\alpha^k>0$ such that
\begin{equation}
    \label{eq:line search 3}
    \begin{aligned}
    \widetilde{F}(\T_{\A}^{k+1})\leq& \widetilde{F}(\T_{\B}^k)+\langle\nabla\widetilde{F}(\T_{\B}^k),\T_{\A}^{k+1}-\T_{\B}^k\rangle\\
    &+\frac{\delta}{2\alpha^k} \|\T_{\A}^{k+1}-\T_{\B}^k\|^2 ,
\end{aligned}
\end{equation}
   it ensures a sufficient decrease of the Lagrangian function~\eqref{eq:new Lagarangian} at $\Sub^k$ along the direction $\T_{\A}^{k+1}-\T_{\B}^k$, that is:  
\begin{equation*}
   \mathcal{L}(\T_{\A}^{k+1},\Sub^k)\leq \mathcal{L}(\T_{\B}^k,\Sub^k)-\left(1-\frac{\delta}{2}\right)\frac{1}{\alpha^k}\|\T_{\A}^{k+1}-\T_{\B}^k\|^2.
\end{equation*}
 Thus, $\T_{\A}^{k+1}-\T_{\B}^k$ becomes our viable primal search direction, and  $\Sub^k$ will serve as the   dual variables. 
 
Using the  $\Sub$-variable and absorbing the   $\Z$-variables,~\eqref{eq:FPI_1} can be rewritten as: given   the largest   $\alpha^k>0$ satisfying~\eqref{eq:line search 3},  \vspace{-0.4cm}
\begin{subequations}\label{eq:ATOS A}
    \begin{align}
\T_{\A}^{k+1}&=\prox_{\alpha^k\widetilde{G}}\left(\T_{\B}^k-\alpha^k\Sub^k-\alpha^k\nabla\widetilde{F}(\T_{\B}^k)\right)
        \label{eq:ATOS X}\\
        \T_{\B}^{k+1}&=\prox_{\alpha^k\widetilde{R}}\left(\T_{\A}^{k+1}+\alpha^k\Sub^k\right),\\
        \label{eq:ATOS S}
\Sub^{k+1}&=\Sub^k+\frac{1}{\alpha^k}(\T_{\A}^{k+1}-\T_{\B}^{k+1}).
    \end{align}
\end{subequations}


\subsection{A decentralized implementation}
The algorithm in~\eqref{eq:ATOS A}   is not fully decentralized, due to the reliance on $\texttt{prox}_{\alpha^k\widetilde{G}}$ and the backtracking procedure \eqref{eq:line search 3} involving the function $\widetilde F$,  which cannot be performed locally and independently by the agents.  In this section we   cope with these issues, and proposed a fully decentralized implementation of the algorithm in~\eqref{eq:ATOS A}.

 Let us partition the variables  as 
\begin{equation}\T_{\A}^k:=\begin{bmatrix}\T_{\A,1}^k\\\T_{\A,2}^k\end{bmatrix},\quad \T_{\B}^k:=\begin{bmatrix}\T_{\B,1}^k\\\T_{\B,2}^k\end{bmatrix},\quad \mathbf{S}^k:=\begin{bmatrix}\mathbf{S}^k_1\\\mathbf{S}^k_2\end{bmatrix}.
\end{equation}  
The update of the $\T_{\B}$-variables involving $\texttt{prox}_{\alpha^k\widetilde{R}}$ reads 
\begin{equation}
    \label{eq:BCV X 1'}
        \T_{\B}^{k+1}
=\begin{bmatrix}
        \texttt{prox}_{\alpha^k\widetilde{R}}(\T_{\A,1}^{k+1}+\alpha^k\Sub^k_1)\\
        \mathbf{0}
    \end{bmatrix}.
\end{equation}
The expression of $\texttt{prox}_{\alpha^k\widetilde{G}}$ is obtained using   the following  

\noindent \textbf{Fact} (e.g.,  \cite[Eq.~(2.6)]{ryu2022large}):  
    \label{lemma:composite G}
Let $g(u)=f^*(A^\top u)$, with  $f:\mathbb{R}^n\to \mathbb{R}\cup\{-\infty, \infty\}$ assumed to be  convex, closed and proper. Suppose,  $\texttt{ri}\,\texttt{dom}\, f^*\cap \texttt{range}(A^\top)\neq\emptyset$. Then $v=\texttt{prox}_{\alpha,g}(u)$ if and only if there exists $x\in \mathbb{R}^n$ such that 
\begin{equation}\label{eq:sub-conj}
       \begin{aligned}
        x&\in\texttt{argmin}_y\{f(y)-\langle u,Ay\rangle+\frac{\alpha}{2}\|Ay\|^2\}\quad \text{and} \\
        v&=u-\alpha Ax.
    \end{aligned}
\end{equation}
Notice that   $\widetilde{G}=(\delta_{\{\mathbf 0\}}^*)^*\circ [\PL,M]$, with $\delta_{\{\mathbf 0\}}^*$ denoting the conjugate of $\delta_{\{\mathbf 0\}}$.  Since  $\delta_{\{\mathbf 0\}}^*$ is convex, closed and proper    and $\texttt{ri}\texttt{dom} \,\delta_{\{\mathbf 0\}}\cap \texttt{Range}([\PL,M]])\neq\emptyset$,  we can   use \eqref{eq:sub-conj} to express $\texttt{prox}_{\alpha^k\widetilde{G}}$. Hence,    \eqref{eq:ATOS X} can be rewritten as   
\begin{equation}
    \label{eq:A}
 \T_{\A}^{k+1}=\begin{bmatrix}
        \T_{\B,1}^k-\alpha^k\Sub_1^k -\alpha\nabla F(\T_{\B,1}^k)\\
        \T_{\B,2}^k-\alpha^k\Sub_2^k  
    \end{bmatrix}-\alpha^k\begin{bmatrix}
        \PL\\ M
    \end{bmatrix}\Y^{k+1},
\end{equation}
where \begin{equation}
    \label{eq:BCV Y 1'}
    \begin{aligned}       \Y^{k+1}\in&\texttt{argmin}_{\Y}\{\delta_{\{\mathbf{0}\}}^*(\Y)-\langle \T_{\B,2}^k-\alpha^k\Sub_2^k,M\Y\rangle\\&-\langle \T_{\B,1}^k-\alpha^k\Sub^k_1-\alpha^k\nabla F(\T_{\B,1}^k),\PL \Y\rangle\\
    &+\frac{\alpha^k}{2}\left(\|\PL \Y\|^2+\|M \Y\|^2\right)\}.
    \end{aligned}
\end{equation}


We proceed choosing $M$ so that  (\ref{eq:BCV Y 1'}) admits a closed form solution, computable locally at the agents' sides.  
Combining~\eqref{eq:ATOS S},~\eqref{eq:BCV X 1'} and~\eqref{eq:A}, it holds
  $\Sub^k_2=-M\Y^{k}.$
Substituting this expression in~\eqref{eq:BCV Y 1'}  and using the fact  
$$
   \delta^*_{\{\mathbf 0\}}(\Y)=\sup_{\Z\in\mathbb{R}^{m\times d}}\left\{\langle\Y,\Z\rangle-\delta_{\{\mathbf{0}\}}(\Z)\right\}=\langle\Y,\mathbf{0}\rangle=0,
 $$ yields 
\begin{equation*}
 \begin{aligned}       &\Y^{k+1}\in
   \texttt{argmin}_\Y\,\left\{ \alpha^k \|\Y-\Y^k\|^2_{\PL^2+M^2}\right.\\
 &-2\langle\PL(\T_{\B,1}^k-\alpha^k\Sub_1^k-\alpha^{k}\PL\Y^k-\alpha^k\nabla F(\T_{\B,1}^k)),\Y\rangle\}.
 \end{aligned}
\end{equation*}
Choosing $M$ to linearize the quadratic term in the objective function above--$M:=\sqrt{I-\PL^2}$--yields 
\begin{equation}
\label{eq:BCV Y 3'}
\begin{aligned}
\Y^{k+1}\!\!=& \Y^k+\frac{1}{\alpha^k}\PL(\T_{\B,1}^k\!-\!\alpha^k\Sub_1^k\!-\!\alpha^{k}\PL\Y^k-\alpha^k\nabla F(\T_{\B,1}^k)).
\end{aligned}
\end{equation}
Introducing the variable  $\D^k:=\PL \Y^k$, and setting  $\PL^2=I-W$, with $W=(1-c)I+c\widetilde{W}\in W_{\mathcal{G}}$ for some $c\in(0,1/2)$, 
\eqref{eq:BCV Y 3'} can be rewritten in terms of the  $\D$-variables as
\begin{equation}
\label{eq:D}
    \D^{k+1}=\D^k+\frac{1}{\alpha^k}(I-W)(\T_{\B,1}^k-\alpha^k(\Sub_1^k+\D^k+\nabla F(\T_{\B,1}^k))),
\end{equation}
and so the update of the    
 $\T_{\A}$-variables (using \eqref{eq:A}):  
\begin{equation}\label{eq:TA-final}
    \T_{\A,1}^k=W\T_{\B,1}^k+W(\nabla F(\T_{\B,1}^k)+\Sub_1^k+\D^k)).
\end{equation}
In summary, the proposed  algorithm is given by   (\ref{eq:BCV X 1'}),  (\ref{eq:D}), and (\ref{eq:TA-final}).  The algorithm is now fully decentralized in the updates of $\T$-and $\D$-variables.   

We are left to address the implementability of the line-search procedure~\eqref{eq:line search 3} over the network.   To this end, we introduce a local stepsize $\alpha_i^k$ for each agent $i$. This stepsize is determined by performing backtracking on the local loss function $f_i$. Specifically, for $t_{\A,i}^{k+1}$ and $t_{\B,i}^k$, which are the $i$-th rows of $\T_{\A,1}^{k+1}$ and $\T_{\B,1}^k$ respectively, each agent finds the largest $\alpha_i^k$ that satisfies: 
\begin{equation}
\begin{aligned}
    \label{eq:line search}
    f_i(t_{\A,i}^{k+1})\leq& f_i(t_{\B,i}^{k})+\langle \nabla f_i(t_{\B,i}^{k}),t_{\A,i}^{k+1}-t_{\B,i}^{k}\rangle\\
    &+\frac{\delta}{2 {\alpha}_i^k}\|t_{\A,i}^{k+1}-t_{\B,i}^{k}\|^2.
\end{aligned}
\end{equation}
Setting $\alpha^k:=\min_{i\in[m]}{\alpha}_i^k$
ensures that the global line-search condition \eqref{eq:line search 3} is also met.  
 
The proposed decentralized algorithm is outlined in Algorithm 1, with the detailed backtracking line-search procedure provided in Algorithm 2. For clarity in presentation, we refer to $\T_{\B,1}^k$ and $\Sub_1^k$ as $\X^k$ and $\Sub^k$, respectively.\smallskip 

\noindent \textit{On the global-min consensus:} Step \texttt{(S.3)} entails executing a min-consensus operation across the network to synchronize a common stepsize, $\alpha^k = \min_{i \in [m]} \alpha_{i}^k$, among the agents. 
This min-consensus protocol adapts seamlessly to contemporary wireless mesh network technologies. In fact, current network infrastructures support diverse communication interfaces, in particular  WiFi and LoRa (Low-Range) \cite{Kim_Lim_Kim_2016, Janssen_BniLam_Aernouts_Berkvens_Weyn_2020}. WiFi provides high-speed, short-range communication capabilities whereas LoRa supports long-range communication over low rates, perfectly suited for broadcasting minimal information to {\it all nodes} across the network in a single hop. In such multi-interface environments, Algorithm 1 is perfectly implementable: it employs WiFi for transmitting vector variables in Steps \texttt{(S.1)} and uses LoRa for executing the min-consensus in Step \texttt{(S.3)}. Additionally, before transmission, the $\alpha_{i}^k$ values can be quantized to lower approximations using  a few bits.

Local-min operations leveraging only neighboring communications will be  considered in Sec.~\ref{sec:local-min}.


\begin{algorithm}[ht!]
\centering
\resizebox{\columnwidth}{!}{%
  \begin{minipage}{\columnwidth}
\caption{Decentralized Adaptive Three Operator Splitting  (global\_DATOS)}
  \noindent \textbf{Data:} (i) Initialization:  $\alpha^{-1} \in (0, \infty)$, $\X^0 \in \mathbb{R}^{m\times n}$, $\Sub^0 \in \mathbb{R}^{m\times n}$  and  $\D^0 = 0$;  (ii) Backtracking   parameters: $\delta \in (0,1]$;  (iii) Gossip matrix  $W=(1-c)I+c\widetilde{W}$, $c\in(0,1/2)$.		
		\begin{algorithmic}[1]
           \State \texttt{(S.1) Communication Step: } $$\X^{k+1/2}=W\X^k,\quad \D^{k+1/2}=W(\nabla F(\X^k)+\Sub^k+\D^k);$$
          
			\State \texttt{(S.2) Decentralized line-search: }Each agent updates $\overline{\alpha}_i^k$ according to:
   	$$\overline{\alpha}_i^k=\texttt{Linesearch}(\alpha^{k-1}, f_i, x_i^k, x_i^{k+1/2},-d^{k+1/2}_i,\gamma_k,\delta);$$
			\State \texttt{(S.3)} \texttt{Global min-consensus: }$$\alpha^k=\min_{i\in[m]}\overline{\alpha}_i^k;$$
     \State \texttt{(S.4) Updates of the primal and dual variables:}
     \begin{align*}
    \X^{k+1}&=\prox_{\alpha^k R}(\X^{k+1/2}-\alpha^k\D^{k+1/2}+\alpha^k\Sub^k)\\
    \Sub^{k+1}&=\Sub^k+\frac{1}{\alpha^k}(\X^{k+1/2}-\X^{k+1}-\alpha^k\D^{k+1/2})\\       \D^{k+1}&=\D^{k+1/2}+\frac{1}{\alpha^k}(\X^k-\X^{k+1/2}-\alpha^k\nabla F(\X^k)-\alpha^{k}\Sub^k).
     \end{align*}
			\State \texttt{(S.5)}    If a termination criterion is not met,  $k\leftarrow k+1$ and go to step \texttt{(S.1)}.
\end{algorithmic}
 \end{minipage}
}
\label{alg:DATOS}
\end{algorithm}
\begin{algorithm}\centering
\resizebox{\columnwidth}{!}{%
  \begin{minipage}{\columnwidth}
\caption{\texttt{Linesearch}($\alpha,f,x_1,x_2,d,\delta$)}
		\begin{algorithmic}[1]
        \State $\alpha^+=\alpha;$
        \State $x^+=x_2+\alpha^+ d$; Set $t=1$;
           \While{$f(x^+)> f(x_1)+\langle\nabla f(x_1),x^+-x_1\rangle+\frac{\delta}{2\alpha^+}\|x^+-x_1\|^2$}
			\State $\alpha^+\leftarrow \alpha^+/2$;
            \State $x^+\leftarrow x_2+\alpha^+ d$;
            \State $t\leftarrow t+1$;
            \EndWhile
            \Return $\alpha^+$.
\end{algorithmic}\label{alg:backtracking_DATOS}
\end{minipage}}
\end{algorithm}
\subsection{Convergence guarantees}
 This section study  convergence of Algorithm~1. For analytical convenience, we will   refer  to the updates written as in  \eqref{eq:ATOS A}. In fact, since   $\D^0=\mathbf{0}\in\texttt{span}(I-W)$, Algorithm~\ref{alg:DATOS} is equivalent to \eqref{eq:ATOS A}, provided that the stepsize $\alpha^k$ therein is updated according to the backtracking procedure described in \texttt{(S.3)} of Algorithm~1. Henceforth, all references to \eqref{eq:ATOS A} will implicitly assume this condition.

Given the   reformulation of~\eqref{eq:ausiliary} based on the Lagrangian function $\mathcal L$ in~\eqref{eq:new Lagarangian} and strong duality (guaranteed  by Assumption~\ref{ass:function}), a saddle point $(\T^\star,\mathbf{S}^\star)$ of $\mathcal L$ is such that $\T^\star$ 
is a minimizer of~\eqref{eq:ausiliary}. For any given saddle-point    $(\T^*,\Sub^*)$  of   $\mathcal L$, let  $\mathcal{B}_T\times \mathcal{B}_S$ be a bounded set containing $(\T^*,\Sub^*)$. We introduce the partial primal-dual gap function~\cite{chambolle2011first},
\begin{equation}
    \label{eq:new lyapunov}
    \mathcal{G}_{\mathcal{B}_T\times \mathcal{B}_S}(\T_\A,\Sub) :=\max_{\Sub\in B_S}\mathcal{L}(\T_{\A},\Sub)-\min_{\T\in\mathcal{B}_T}\mathcal{L}(\T,\Sub).
\end{equation} 
Notice that   $\mathcal{G}_{\mathcal{B}_T\times \mathcal{B}_S}(\T_\A,\Sub)\geq 0$, for any  $(\T_\A,\Sub)\in \texttt{dom}(\mathcal L)$. Yet, $\mathcal{G}_{\mathcal{B}_T\times \mathcal{B}_S}$ is not a valid measure of optimality because   it might vanish for $(\T_\A,\Sub)$ not being a saddle-point of $\mathcal L$ \cite{chambolle2011first}.  However,   if $\mathcal{G}_{\mathcal{B}_T\times \mathcal{B}_S}(\T_\A,\Sub)=0$ for   $(\T_\A,\Sub)\in \texttt{dom}(\mathcal L)$  and, in addition, $(\T_\A,\Sub)$ lies in the interior of $\mathcal{B}_T\times \mathcal{B}_S$, then $(\T_\A,\Sub)$ is a saddle-point of $\mathcal L$. 
The following theorem   shows that the iterates $(\T_\A^k,\Sub^k)\in \texttt{dom}(\mathcal L)$   remain bounded. Hence, one can choose a ``sufficiently large''  $\mathcal{B}_T\times \mathcal{B}_S$ such that $\mathcal{G}_{\mathcal{B}_T\times \mathcal{B}_S}(\T_\A^k,\Sub^k)$ is a valid measure of 
optimality for $(\T_\A^k,\Sub^k)$. 
 We provide next the convergence rate decay of $\mathcal{G}_{\mathcal{B}_T\times \mathcal{B}_S}$ along the    ergodic sequences   of the   iterates   $\{\T_{\A}^k\}$ and  $\{\Sub^k\}$, defined as   $$\!\!\!\overline{\T}^k_{\A}:=\frac{1}{\theta^k}\sum_{t=1}^{k}\alpha^{t-1}\T_{\A}^{t},\,\, \overline{\Sub}^{k}:=\frac{1}{\theta^k}\sum_{t=1}^{k}\alpha^{t-1}\Sub^{t}, \,\,\theta^k:=\sum_{t=1}^{k}\alpha^{t-1}.$$ 
\begin{theorem}
\label{thm:asymptotic convergence}
  Under Assumptions~\ref{ass:function}-\ref{ass:W}, the following holds along the iterates of Algorithm~\ref{alg:DATOS}:
  
  \textbf{(a)} The sequence  $\{\T_{\A}^k,\Sub^k\}_k$ is bounded; and 
  
  \textbf{(b)} For any given saddle-point $(\T^*,\Sub^*)$ of $\mathcal L$, there exists a   bounded set   $\mathcal{B}_{T}\times\mathcal{B}_{S}$ [see \eqref{eq:set_Bs} in the proof for the explicit expression] such that (i)  $(\T^*,\Sub^*)\in \mathcal{B}_{T}\times\mathcal{B}_{S}$, (ii) $(\T_{\A}^k,\Sub^k)\in \texttt{int}(\mathcal{B}_{T}\times\mathcal{B}_{S})$, for all $k\geq 0$; and   (iii)
  \begin{equation}
        \label{eq:sublinea DATOS}
\mathcal{G}_{\mathcal{B}_T\times \mathcal{B}_S}(\overline{\T}_{\A}^k,\overline{\Sub}^k)\leq\frac{1}{k}\cdot\frac{8R_0^2}{ \min\left(\alpha^{-1},\delta/(2{L}^k)\right)},\end{equation}
where ${L}^k$ is the Lipschitz constant of $\nabla F$ over the convex hull of  $\cup_{t=1}^k[\T_{\B,1}^{t-1},\T_{\A,1}^t]$, and  \begin{equation}\label{eq:R-square}R_0^2:=\|\T_{\B}^{0}-\T^*\|^2+(\alpha^{0})^2\|\Sub^{0}-\Sub^*\|^2.\end{equation}

  
\end{theorem}


Theorem~1 establishes the sublinear convergence rate of Algorithm~\ref{alg:DATOS}, matching   $\mathcal{O}(1/k)$  rate of  state-of-art  algorithms solving decentralized convex optimization  problems.  Crucially, Algorithm~\ref{alg:DATOS}  removes the requirement  for prior knowledge of optimization and network parameters, which is common to all prior native decentralized algorithms  applicable to Problem~\eqref{eq:problem}, such as \cite{shi2015proximal,guo2023decentralized}.
Notably   the convergence rate~\eqref{eq:sublinea DATOS} depends on the \textit{local} Lipschitz constants  $\{L^k\}$, defined on the convex hull of the trajectory traveled by the algorithm. This local constant often proves much smaller than the \textit{global} Lipschitz constant $L$,  governing the rates of existing (non adaptive) decentralized algorithms.  This advantage showcases the algorithm's ability to adapt to the local geometry of the optimization function.  
Finally,   the constant $R_0^2$   in (\ref{eq:sublinea DATOS}) is to be expected, it quantifies the ``distance'' of the initial point   from the targeted saddle-point.  

\subsection{Proof sketch of Theorem 1}
\textit{Proof of statement (a):}  It can be shown that   the line-search   terminates in a finite number of iteration, with $\alpha^{k}$ satisfying \begin{equation}\label{eq:lower-bound-step}
\min\left(\alpha^{k-1},\frac{\delta}{2\widetilde L^k}\right)\leq\alpha^k\leq \alpha^{k-1},\end{equation}
where 
$\widetilde L^k\in(0,\infty)$ is the Lipschitz constant of $\nabla F$ over the segment $[\T_{\B,1}^k,\T_{\A,1}^{k+1}]$.    
 Furthermore, 
\begin{equation}
    \label{eq:line search 4}
\begin{aligned}\Delta^k:=&F(\T_{\B,1}^k)+\langle\nabla F(\T_{\B,1}^k),\T_{\A,1}^{k+1}-\T_{\B,1}^k\rangle\\
   &+\frac{\delta}{2\alpha^k}\|\T_{\A,1}^{k+1}-\T_{\B,1}^k\|^2-F(\T_{\A,1}^{k+1})\geq 0.
   \end{aligned}
\end{equation}

 We study the descent of the gap function $\mathcal{G}^{k}$.
    By definition ~\eqref{eq:new lyapunov} of $\mathcal L$,   we have: for any fixed $(\T,\Sub)\in\texttt{dom}(\mathcal{L})$,
    \begin{equation*}
       \begin{aligned}          &\mathcal{L}(\T_{\A}^{k+1},\Sub)-\mathcal{L}(\T,\Sub^{k+1})\\
            \leq&\underbrace{\widetilde{F}(\T_{\A}^{k+1})-\widetilde{F}(\T_{\B}^k)}_{\texttt{term I}}+\underbrace{\widetilde{F}(\T_{\B}^k)-\widetilde{F}(\T)}_{\texttt{term II}}\\
         \end{aligned}\end{equation*}\begin{equation}
            \label{eq:first equation}  \begin{aligned}
            &+\underbrace{\widetilde{G}(\T_{\A}^{k+1})-\widetilde{G}(\T)}_{\texttt{term III}}+\underbrace{\langle\T_{\A}^{k+1}-\T,\Sub^{k+1}\rangle}_{\texttt{term IV}}\\
            &+\underbrace{\widetilde{R}^*(\Sub^{k+1})-\widetilde{R}^*(\Sub)}_{\texttt{term V}}+\underbrace{\langle\T_{\A}^{k+1},\Sub-\Sub^{k+1}\rangle}_{\texttt{term VI}}.
       \end{aligned}
    \end{equation}
    We proceed bounding each term separately.  
    By~\eqref{eq:line search 4}
    \begin{equation}
        \label{eq:term I DATOS}
        \begin{aligned}
            \texttt{term I}
            &\leq
             \langle\nabla F(\T_{\B,1}^k),\T_{\A,1}^{k+1}-\T_{\B,1}^k\rangle\\&\quad +\frac{\delta}{2\alpha^k}\|\T_{\A,1}^{k+1}-\T_{\B,1}^k\|^2.
        \end{aligned}
    \end{equation}
    Using convexity of $F$, yields
    \begin{equation}
        \label{eq:term II DATOS}
       \begin{aligned}
            \texttt{term II}
            \leq&\langle\nabla \widetilde{F}(\T_{\B}^k),\T_{\B}^k-\T\rangle.
       \end{aligned}
    \end{equation}
    From the update  \eqref{eq:ATOS A}, we have
  $$  
        \frac{1}{\alpha^k}(\T_{\B}^k-\alpha^k\Sub^k-\alpha^k\nabla\widetilde{F}(\T_{\B}^k)-\T_{\A}^{k+1})\in\partial \widetilde{G}(\T_{\A}^{k+1}),
    $$
   which, together with  convexity of $\widetilde{G}$, yields
    \begin{equation}
        \label{eq:term III DATOS}
        \begin{aligned}
          \texttt{term III}\leq &\frac{1}{\alpha^k}\langle\T_{\B}^k-\T_{\A}^{k+1},\T_{\A}^{k+1}-\T\rangle\\
          &-\langle\Sub^k+\nabla\widetilde{F}(\T_{\B}^k),\T_{\A}^{k+1}-\T\rangle.  
        \end{aligned}
    \end{equation}
    Using~\eqref{eq:term I DATOS},~\eqref{eq:term II DATOS} and~\eqref{eq:term III DATOS}, we conclude 
    \begin{equation}
        \label{eq:first half}
        \begin{aligned}
            &\texttt{term I}+\texttt{term II}+\texttt{term III}+\texttt{term IV}\\
            \overset{~\eqref{eq:ATOS S}}{\leq}&
            \frac{\delta}{2\alpha^k}\|\T_{\A,1}^{k+1}-\T_{\B,1}^k\|^2
            +\frac{1}{\alpha^k}\langle\T_{\B}^k-\T_{\A}^{k+1},\T_{\A}^{k+1}-\T\rangle\\
            &-\frac{1}{\alpha^k}\langle\T_{\B}^{k+1}-\T_{\A}^{k+1},\T_{\A}^{k+1}-\T\rangle\\
            =&
            \frac{1}{2\alpha^k}\left(\|\T_{\B}^k-\T\|^2-\|\T_{\B}^k-\T_{\A}^{k+1}\|^2+\delta\|\T_{\B,1}^k-\T_{\A,1}^{k+1}\|^2\right.\\
            &\left.-\|\T_{\B}^{k+1}-\T\|^2+\|\T_{\B}^{k+1}-\T_{\A}^{k+1}\|^2\right)\\
            \leq&
            \frac{1}{2\alpha^k}\left(\|\T_{\B}^k-\T\|^2-\|\T_{\B}^{k+1}-\T\|^2+\|\T_{\B}^{k+1}-\T_{\A}^{k+1}\|^2\right),
        \end{aligned}
    \end{equation}
     where in the last inequality we used  $0<\delta< 1$. 
     
     Using~\eqref{eq:ATOS X} and~\eqref{eq:ATOS S}, it follows
    \begin{equation}
        \label{eq:R subgradient}
        \Sub^{k+1}\in\partial\widetilde{R}(\T_{\B}^{k+1})\Rightarrow\T_{\B}^{k+1}\in\partial\widetilde{R}^*(\Sub^{k+1}).
    \end{equation}
    Leveraging convexity of $\widetilde{R}^*$, we have  
    \begin{equation}
        \label{eq:term V DATOS}
\texttt{term V}\leq \langle\T_{\B}^{k+1},\Sub^{k+1}-\Sub\rangle.
    \end{equation}
    Thus, we can bound  $\texttt{term V}+\texttt{term VI}$ as\begin{equation}
        \label{eq:second half}
        \begin{aligned}
            &\texttt{term V}+\texttt{term VI}
            \leq \langle \T_{\B}^{k+1}-\T_{\A}^{k+1},\Sub^{k+1}-\Sub\rangle\\
\overset{\eqref{eq:ATOS S}}{=}&\frac{1}{2\alpha^k}\left(\|\alpha^k(\Sub^k-\Sub)\|^2
-\|\alpha^k(\Sub^{k+1}-\Sub)\|^2-\right.\\&\left.\|\T_{\B}^{k+1}-\T_{\A}^{k+1}\|^2\right).
        \end{aligned}
    \end{equation}
Combining~\eqref{eq:first half} and~\eqref{eq:second half} and using  $\T_{\B,2}^k=\mathbf{0}$, yields
 \begin{equation}
        \label{eq:DATOS: recursive}
        \begin{aligned}
& 2\alpha^k\left(\mathcal{L}(\T_{\A}^{k+1},\Sub)-\mathcal{L}(\T,\Sub^{k+1})\right)\\
 \leq&
 \|\T_{\B}^k-\T\|^2+(\alpha^k)^2\|\Sub^{k}-\Sub\|^2\\
 &-(\|\T_{\B}^{k+1}-\T\|^2+(\alpha^k)^2\|\Sub^{k+1}-\Sub\|^2).
        \end{aligned}
    \end{equation}
 Telescoping~\eqref{eq:DATOS: recursive} and leveraging 
  $\alpha^{k}\leq\alpha^{k-1}$ yield
  \begin{equation}
      \label{eq:telescope bound}
      \begin{aligned}
&\|\T_{\B}^k-\T\|^2+(\alpha^k)^2\|\Sub^k-\Sub\|^2\\
&+2\sum_{t=1}^{k}\alpha^t\left(\mathcal{L}(\T_{\A}^{t},\Sub)-\mathcal{L}(\T,\Sub^{t})\right)\\ \leq& \|\T_{\B}^0-\T\|^2+(\alpha^{0})^2\|\Sub^0-\Sub\|^2.\end{aligned} 
  \end{equation}
  
 Let    
$(\T^*,\Sub^*)$  be  any saddle-point of $\mathcal L$;  set  in (\ref{eq:telescope bound}) $(\T,\Sub)=(\T^*,\Sub^*)$.  We have  $\mathcal{L}(\T_{\A}^{k+1},\Sub^\star)-\mathcal{L}(\T^\star,\Sub^{k+1})\geq 0$.  
     Invoking the fact that $\alpha^k\geq 0$ for any $k$,  
     we conclude 
     \begin{equation}
         \label{eq:T_B bound}
      \begin{aligned}
&\|\T_{\B}^k-\T^\star\|^2+(\alpha^k)^2\|\Sub^k-\Sub^\star\|^2\\
& \qquad  \leq \|\T_{\B}^0-\T^\star\|^2+(\alpha^{0})^2\|\Sub^0-\Sub^\star\|^2,\quad \forall k\geq 0.\end{aligned} 
  \end{equation}
     Hence, $\{\T_{\B}^{k}, \Sub^k\}_k$ is bounded. And so is   $\{\T_{\A}^k\}_k$, due to ~\eqref{eq:ATOS S}. Specifically, given   $R_0^2>0,$    defined  in~\eqref{eq:R-square}, we have
  \begin{equation}
      \label{eq:T_A bound}
      \|\T_{\A}^k-\T^*\|^2+(\alpha^0)^2\|\Sub^k-\Sub^*\|^2\leq 4R_0^2,\quad \forall k\geq 0.
  \end{equation}


           
           \smallskip

          \textit{Proof of statement (b):}
From~\eqref{eq:T_B bound}-\eqref{eq:T_A bound}  it follows that we can define the bounded set $\mathcal{B}_{T}\times\mathcal{B}_{S}$ containing $(\T^*,\Sub^*)$ as
\begin{equation}\label{eq:set_Bs}\mathcal{B}_T\times \mathcal{B}_S = \left\{(\T,\Sub)|\|\T-\T^*\|^2+(\alpha^{0})^2\|\Sub-\Sub^*\|^2\leq 4R_0^2\right\}.\end{equation}
This also ensures $(\T_{\A}^k,\Sub^k)\!\in\! \texttt{int}(\mathcal{B}_{T}\times\mathcal{B}_{S})$,  for any $k\geq 0$. 

Dividing   both   sides of~\eqref{eq:telescope bound} by $2\theta_k$ and leveraging convexity of $\mathcal{L}(\bullet,\Sub)-\mathcal{L}(\T,\bullet)$, yields 
\begin{equation}
    \label{eq:jensen}
 \mathcal{L}(\overline{\T}_{\A}^k,\Sub)-\mathcal{L}(\T,\overline{\Sub}^k)\leq\frac{1}{2\theta^k}\left(\|\T_{\B}^0-\T\|^2+(\alpha^0)\|\Sub^0-\Sub\|^2\right).
\end{equation}
Taking the maximum of $(\T,\Sub)$ over $\mathcal{B}_T\times\mathcal{B}_S$ on both sides of~\eqref{eq:jensen} and using~\eqref{eq:new lyapunov}, we infer
    \begin{equation*}
       \begin{aligned}    
\mathcal{G}_{\mathcal{B}_T\times \mathcal{B}_S}(\overline{\T}_{\A}^k,\overline{\Sub}^k)\leq&
       \frac{D(\mathcal{B}_{T},\mathcal{B}_{S})}{2\theta^k},
       \end{aligned}
    \end{equation*}
 where $D(\mathcal{B}_{T},\mathcal{B}_{S}):=\!\!\sup_{(\T,\Sub)\in\mathcal{B}_{T}\times\mathcal{B}_{S}}\|\T_{\B}^0-\T\|^2+(\alpha^{0})^2\|\Sub^0-\Sub\|^2.$ Moreover, 
 \begin{equation*}
     \theta^k\geq \sum_{t=1}^k\min\left(\alpha^{t-2},\delta/2\widetilde L^{t-1}\right)\geq\min\left(k\alpha^{-1},k\delta/(2{L}^k)\right).
 \end{equation*}   
 Finally, it is not difficult to check that ${D(\mathcal{B}_{T},\mathcal{B}_{S})}=8R^2_0.$ 
 This completes the proof. $\hfill\square$%

\section{From global to local  min-consensus}\label{sec:local-min}
In this section we introduce a variant of Algorithm~\ref{alg:DATOS} wherein the global min-consensus step   \texttt{(S.3)} is replaced by a local one--see Algorithm~3.  
Step \texttt{(S.3)} therein requires now only {\it local} communications with neighboring nodes. 
The update  \texttt{(S.3)} produces now different local stepsizes, collected in    $ \Lambda^k:=\texttt{diag}(\alpha_1^k,\alpha_2^k,\cdots,\alpha_m^k)$. 
 To simplify the notation,   we defined    $\prox_{\Lambda^k,R}(\X):=
    [\prox_{\alpha_1^k r}(x_1),$ $
    \cdots,
    \prox_{\alpha_m^k r}(x_m)]^{\top}.$

\begin{algorithm}[t!]
\centering
\resizebox{\columnwidth}{!}{%
  \begin{minipage}{\columnwidth}
\caption{Decentralized Adaptive Three Operator Splitting with local min-consensus (local\_DATOS)}
  \noindent \textbf{Data:} (i) initialization:   $\alpha^{-1} \in (0, \infty)$, $\X^0 \in \mathbb{R}^{m\times n}$, $\Sub^0 \in \mathbb{R}^{m\times n}$  and  $\D^0 = 0$;  (ii) Backtracking   parameters $\delta \in (0,1]$; (iii) Gossip matrix $W=(1-c)I+c\widetilde{W}$, $c\in(0,1/2)$.		
		\begin{algorithmic}[1]
           \State \texttt{(S.1) Communication Step: } $$\X^{k+1/2}=W\X^k,\quad \D^{k+1/2}=W(\nabla F(\X^k)+\Sub^k+\D^k);$$
          
			\State \texttt{(S.2) Decentralized line-search: }Each agent updates $\alpha_i^k$ according to:
   	$$\overline{\alpha}_i^k=\texttt{Linesearch}(\alpha_i^{k-1}, f_i, x_i^k, x_i^{k+1/2},-d^{k+1/2}_i,\gamma_k,\delta);$$

			\State \texttt{(S.3)} \texttt{Local min-consensus: }Each agent updates $\alpha_i^k$ according to:
            $$\alpha_i^k=\min_{j\in\mathcal{N}_i}\overline{\alpha}_j^k;$$
            \State\texttt{(S.4) Extra scalar communication step} Let $\Lambda^k:=\texttt{diag}(\alpha_1^k,\alpha_2^k,\cdots,\alpha^k_m)$,
            $$\D_{\Lambda}^k=(I-W)(\Lambda^k)^{-1}\X^k;$$
             
     \State \texttt{(S.5) Updates of the auxiliary, dual and primal variables:}
     \begin{align*}
    \X^{k+1}&=\prox_{\Lambda^k, R}(\X^{k+1/2}-\Lambda^k\D^{k+1/2}+\Lambda^k\Sub^k)\\
    \Sub^{k+1}&=\Sub^k+(\Lambda^k)^{-1}(\X^{k+1/2}-\X^{k+1}-\D^{k+1/2})\\
\D^{k+1}&=\D^{k+1/2}+\D_{\Lambda}^k-\nabla F(\X^k)-\Sub^k.
     \end{align*}
			\State \texttt{(S.6)}    If a termination criterion is not met,  $k\leftarrow k+1$ and go to step \texttt{(S.1)}.
\end{algorithmic}
\end{minipage}
}
\label{alg:DATOS_local}
\end{algorithm}

Notice that  Algorithm~3 differs from Algorithm 2 also in the updates of the $\D$-variables, due to the fact that  $\Lambda^k\notin \texttt{span}(I)$, for all $k$.   In fact, in such a setting, one can show that the update of the  $\D$-variables should  be instead  
\begin{equation}
    \label{eq:BCV D 3''}
\D^{k+1}\!\!=\! \D^k+\PL(\Lambda^k)^{-1}\PL(\X^k-\Lambda^{k}\Sub^k-\Lambda^{k}\D^k-\Lambda^k\nabla F(\X^k)).
\end{equation} However, this update is not implementable on a network because $\PL(\Lambda^k)^{-1}\PL$ is not compliant with the graph $\mathcal G$. The proposed approach is then to ``approximate'' \eqref{eq:BCV D 3''} by 
\begin{equation}
    \label{eq:BCV D 3'''}
    \D^{k+1}\!\!=\! \D^k+\PL^2(\Lambda^k)^{-1}(\X^k-\Lambda^{k}\Sub^k-\Lambda^{k}\D^k-\Lambda^k\nabla F(\X^k)).
\end{equation}
Setting $\PL^2=I-W$,   \eqref{eq:BCV D 3'''} is computable on the network. The update \eqref{eq:BCV D 3'''} serves as an effective approximation of \eqref{eq:BCV D 3''} in the sense that it ensures   two critical properties for convergence:  \textbf{i)}   $\D^k\in\texttt{span}(I-W)$, for all $k$; and \textbf{ii)}   if $\Lambda^k\in\texttt{span}(I)$, (\ref{eq:BCV D 3'''}) recovers (\ref{eq:BCV D 3''}). A key property of the local-min consensus coupled with the proposed line-search procedure is    to ensure that $\Lambda^k$ will fall within $\texttt{span}(I)$ after a finite number of iterations. 



Convergence of Algorithm~\ref{alg:DATOS_local} is summarized    in Theorem~\ref{thm:convergence local_DATOS}, whose proof  is omitted because of space limit--see~\cite{chen-et_al-SIOPT25}.   Similarly to Theorems 1, it is convenient to state Theorem 2 referring to the equivalent formulation of Algorithm~3,  using the  updates in terms of  $\T_{\A,1}^k$ and $\T_{\B,2}^k$, as given in~\eqref{eq:TA-final} and~\eqref{eq:BCV X 1'}, respectively,  with the stepsize $\alpha^k I$ therein  replaced by  $\Lambda^k$, and updated according to (\texttt{S.2})-(\texttt{S.3}) in Algorithm 3. In Theorem~\ref{thm:convergence local_DATOS}, the  $D(\mathcal{B}_{T},\mathcal{B}_{S})$ is defined as in    Theorem~\ref{thm:asymptotic convergence}.

\begin{figure*}[htbp] 
    \centering
    \includegraphics[width=0.3\linewidth]{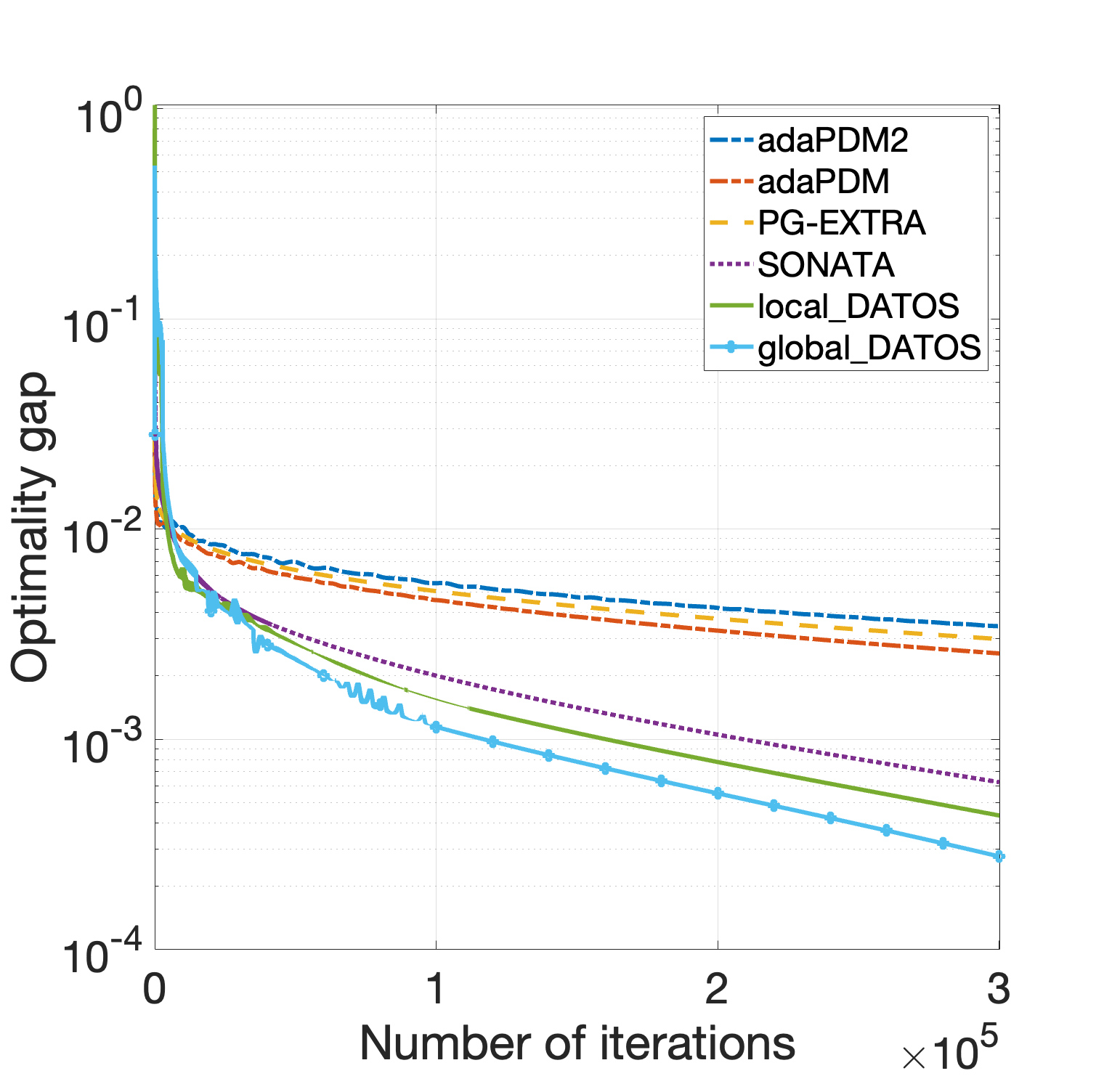}
    \hfill
    \includegraphics[width=0.3\linewidth]{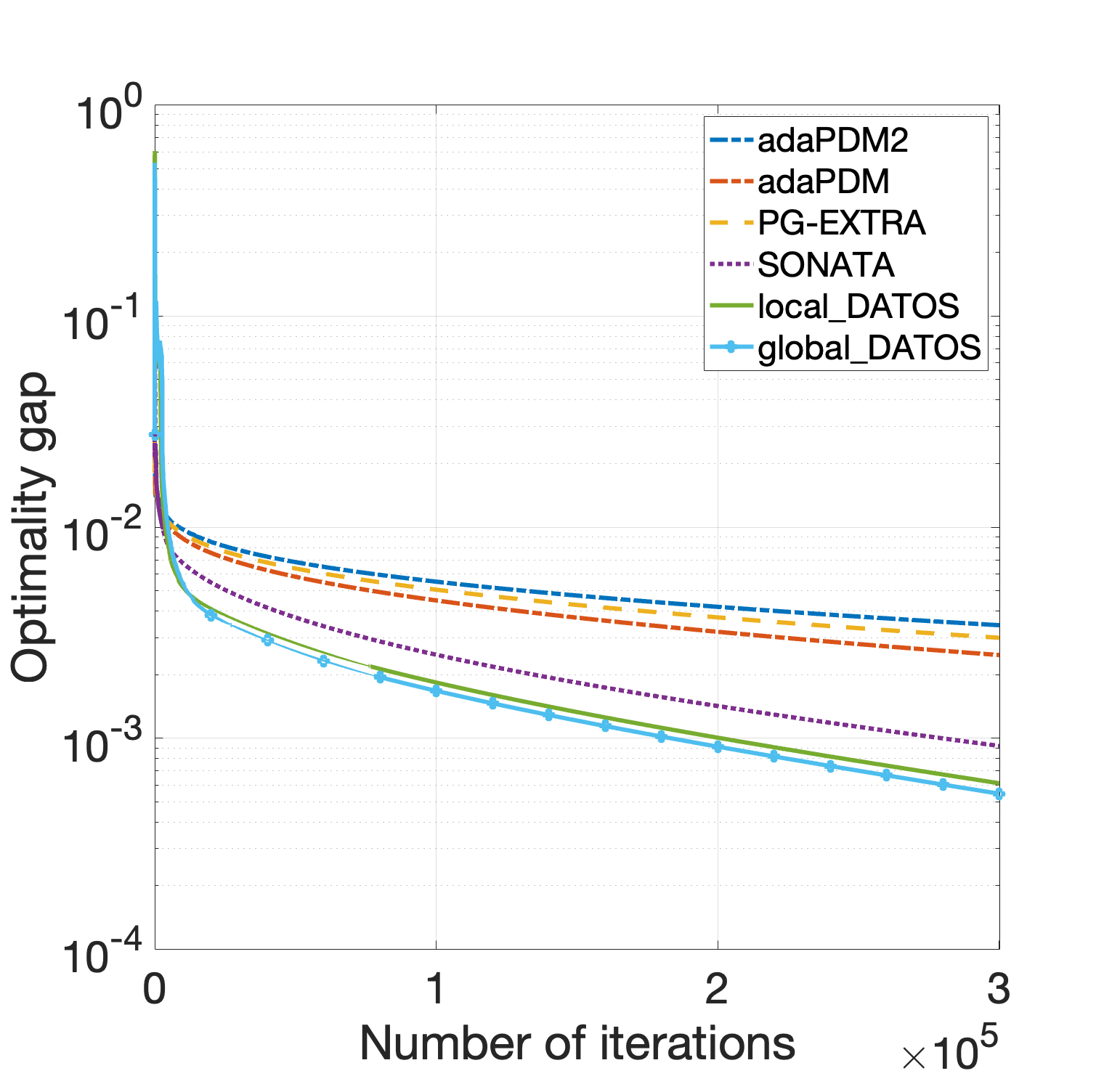}
    \hfill
    \includegraphics[width=0.3\linewidth]{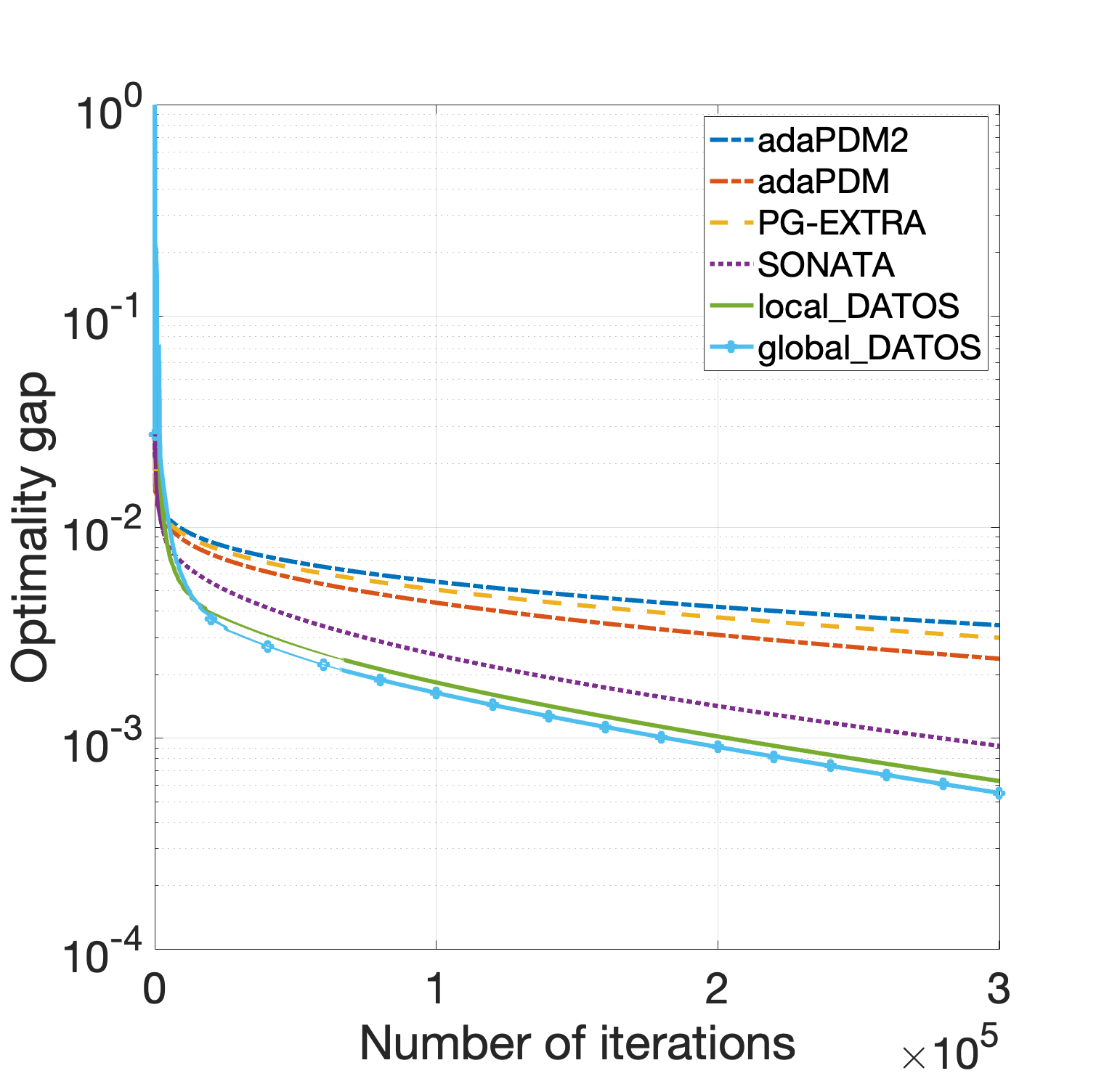}
\caption{Logistic regression with $\ell_1$-regularization: ${\frac{1}{m}\sum_{i=1}^m u(x_i)-u(x^*)}$ v.s. \# iterations.  Comparison of PG-EXTRA, SONATA, adaPDM, adaPDM2, global\_DATOS and local\_DATOS on Erdos-Renyi graphs with different edge-probability:  $p=0.1$  (left);   $p=0.5$ (middle); and  $p=0.9$ (right).}\vspace{-0.45cm} 
\label{fig:general cvx}
\end{figure*}

\begin{figure*}[htbp] 
    \centering
    \includegraphics[width=0.3\linewidth]{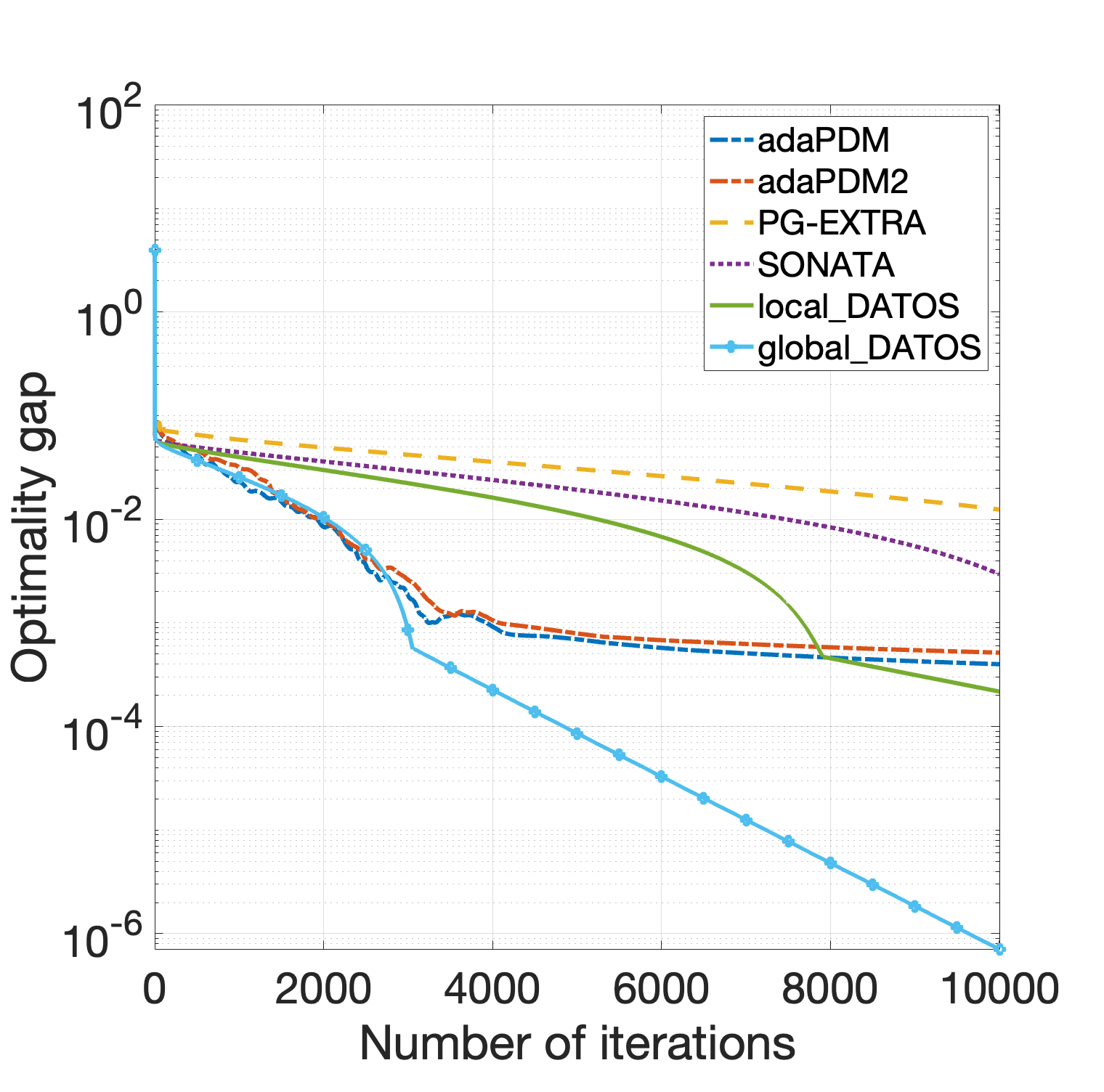}
    \hfill
    \includegraphics[width=0.3\linewidth]{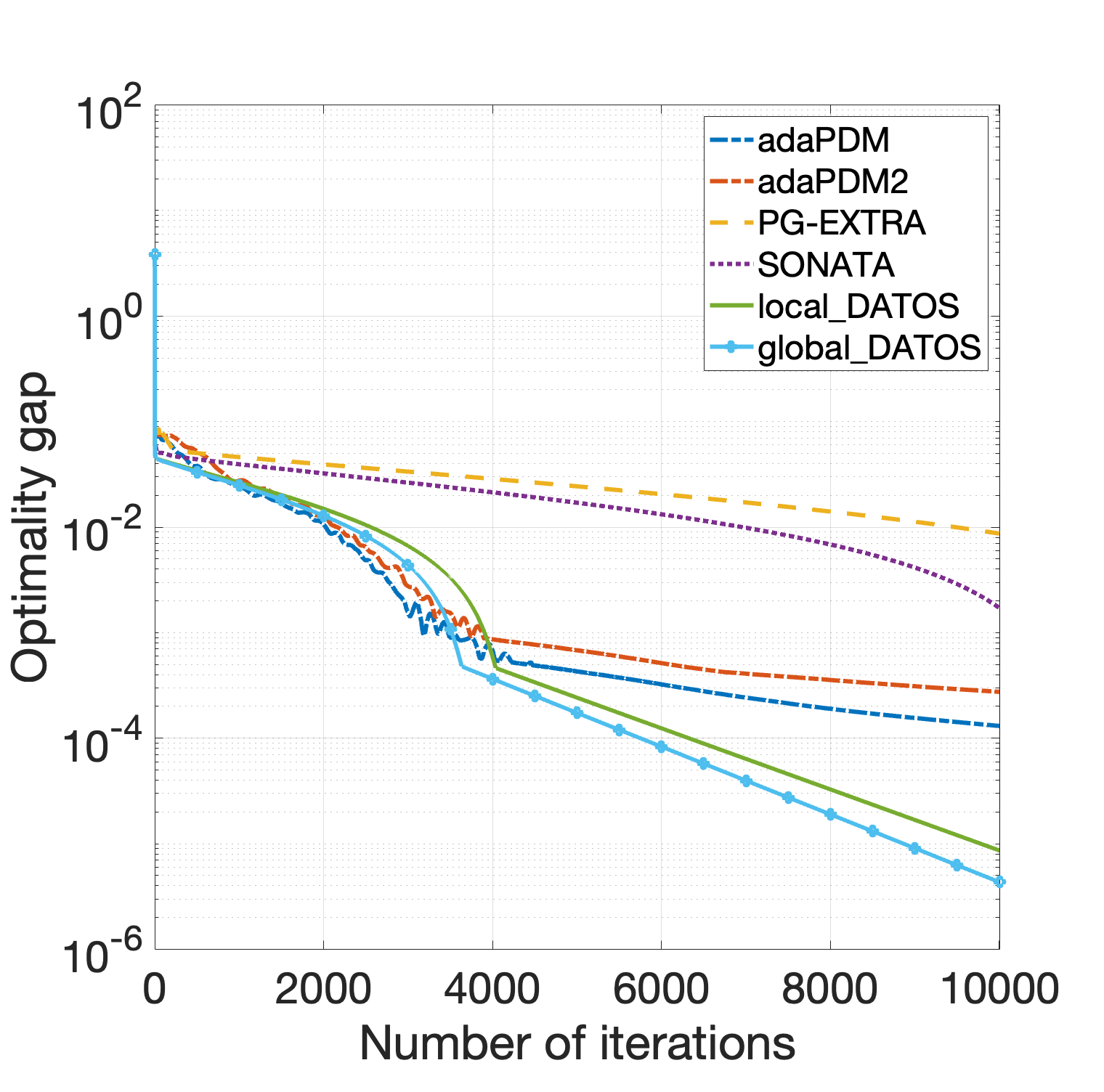}
    \hfill
    \includegraphics[width=0.3\linewidth]{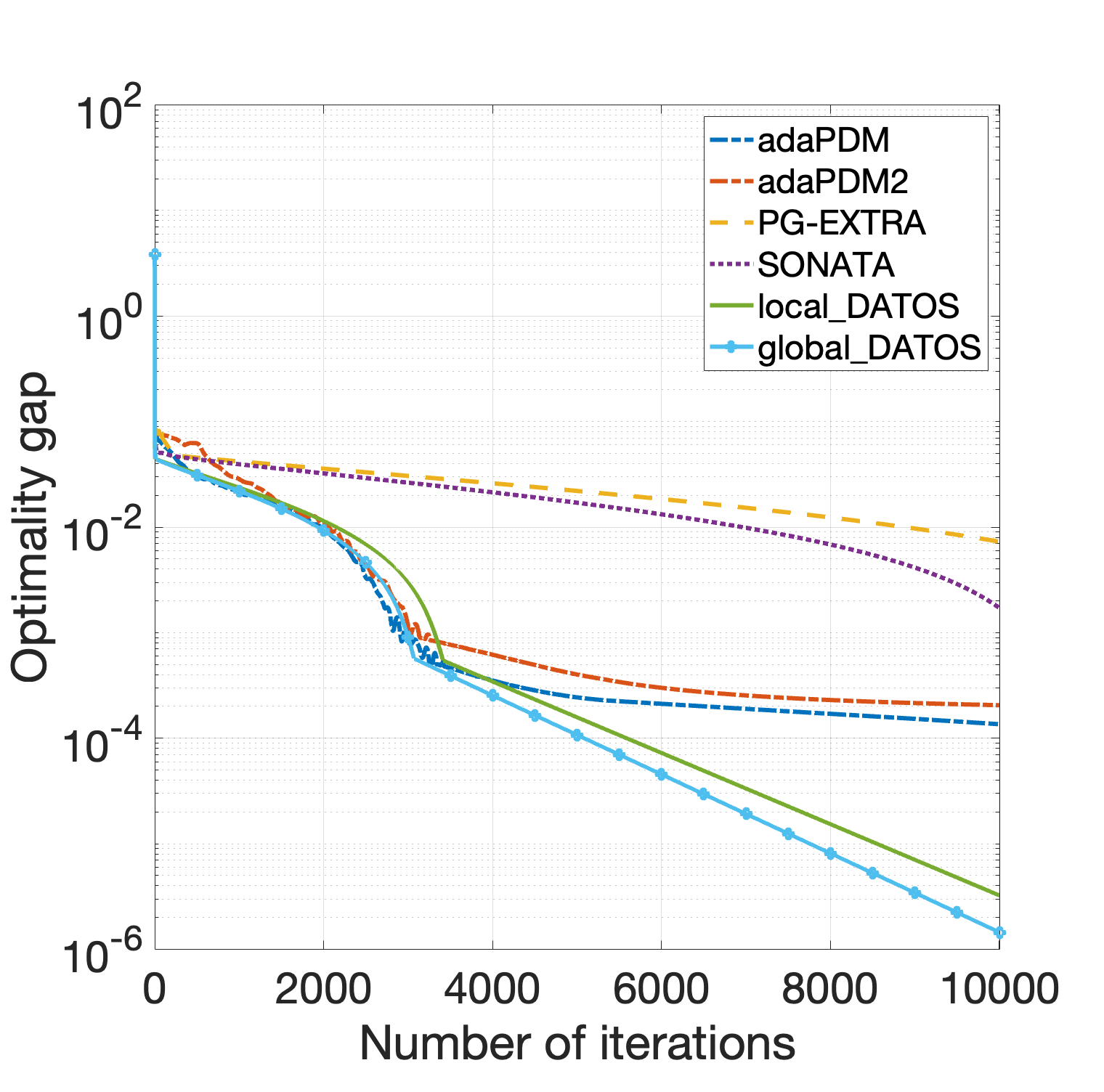}
    \caption{Maximum likelihood estimate of the covariance matrix: ${\frac{1}{m}\sum_{i=1}^m u(x_i)-u(x^*)}$ v.s. \# iterations.  Comparison of PG-EXTRA, SONATA, adaPDM, adaPDM2, global\_DATOS and local\_DATOS on Erdos-Renyi graphs with different edge-probability:  $p=0.1$  (left);   $p=0.5$ (middle); and  $p=0.9$ (right).}\vspace{-0.45cm}

\label{fig:general local smooth}
\end{figure*}

\begin{figure*}[htbp]
    \centering
    \includegraphics[width=0.3\linewidth]{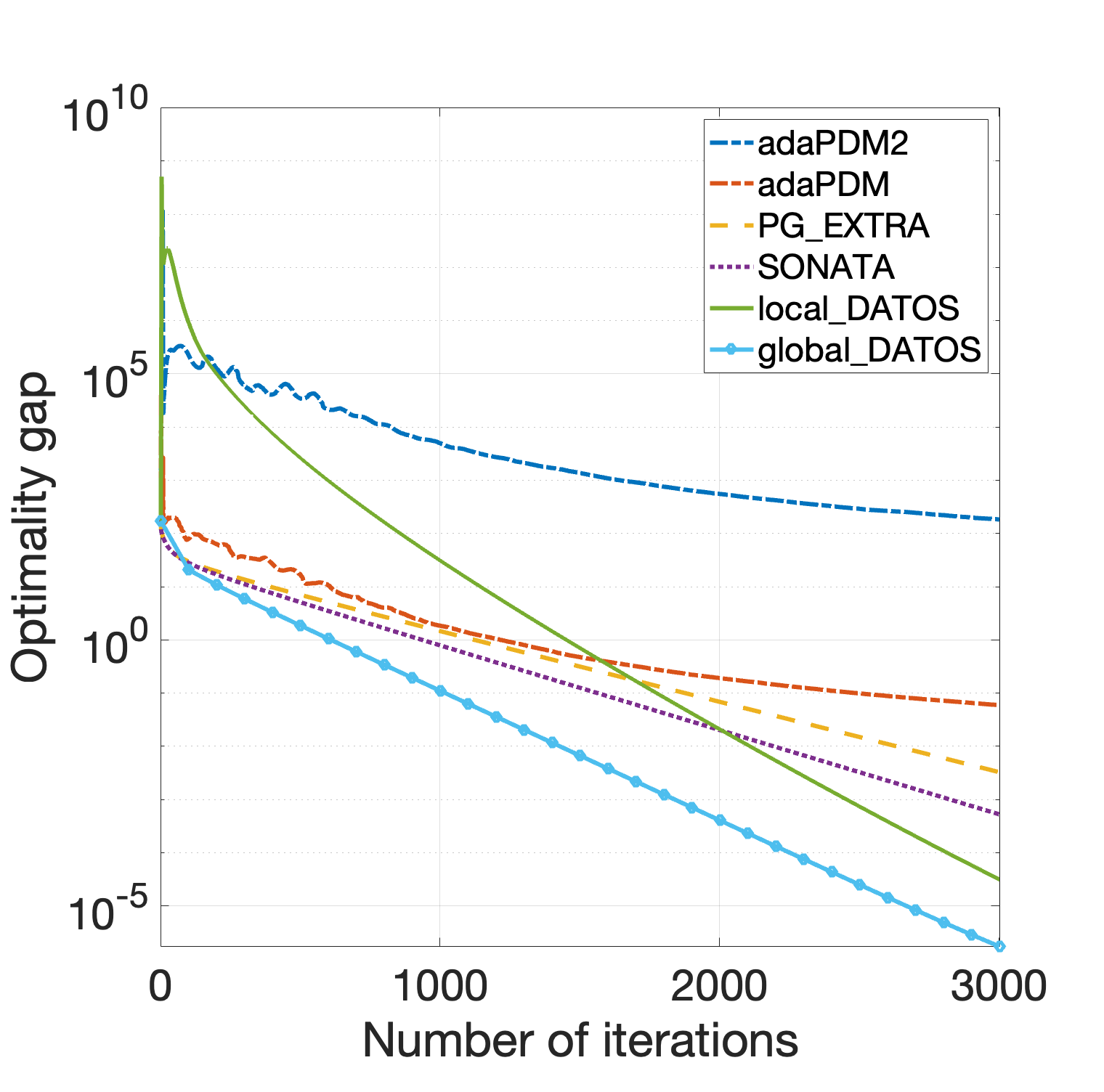}
    \hfill
    \includegraphics[width=0.3\linewidth]{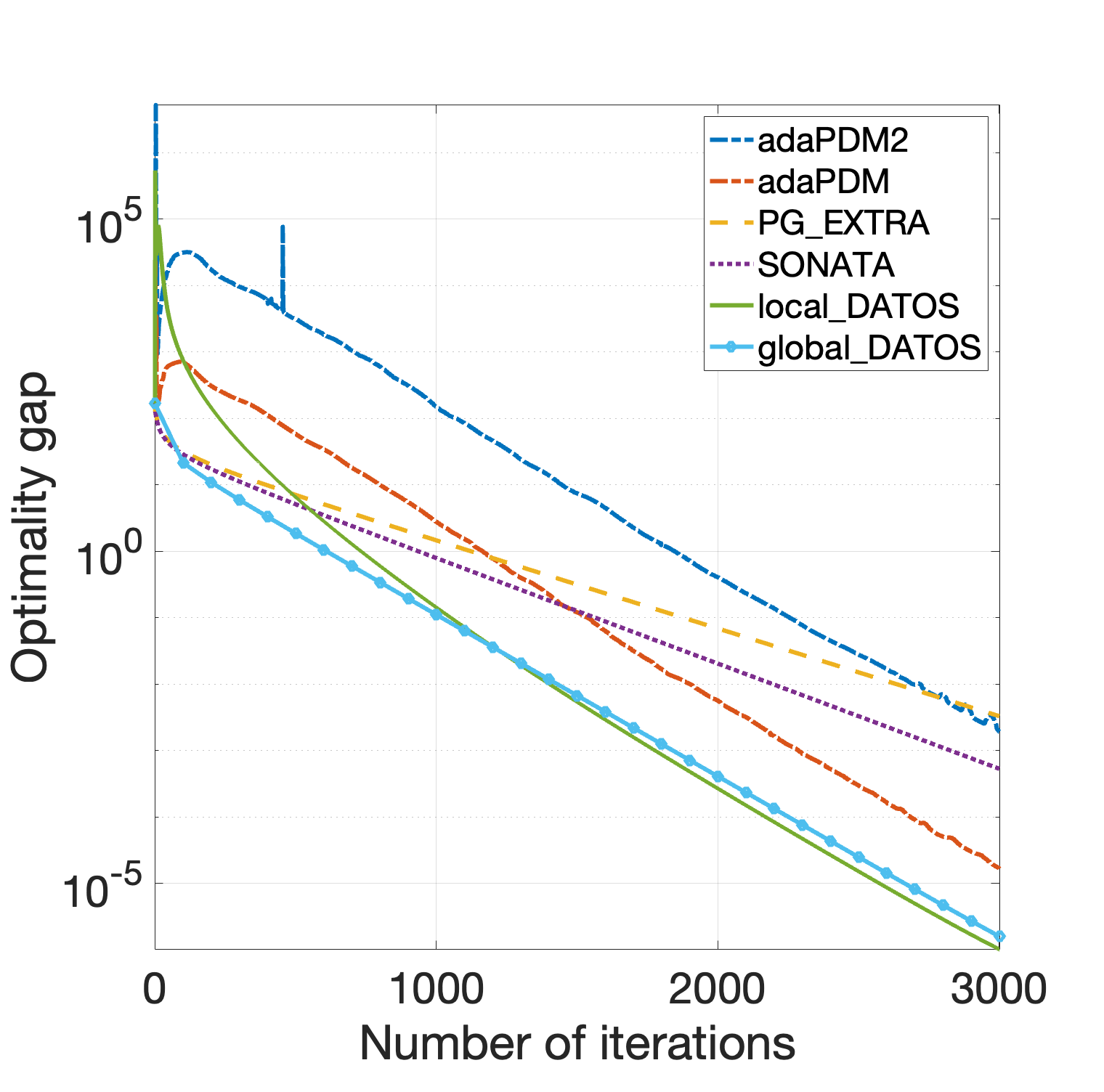}
    \hfill
    \includegraphics[width=0.3\linewidth]{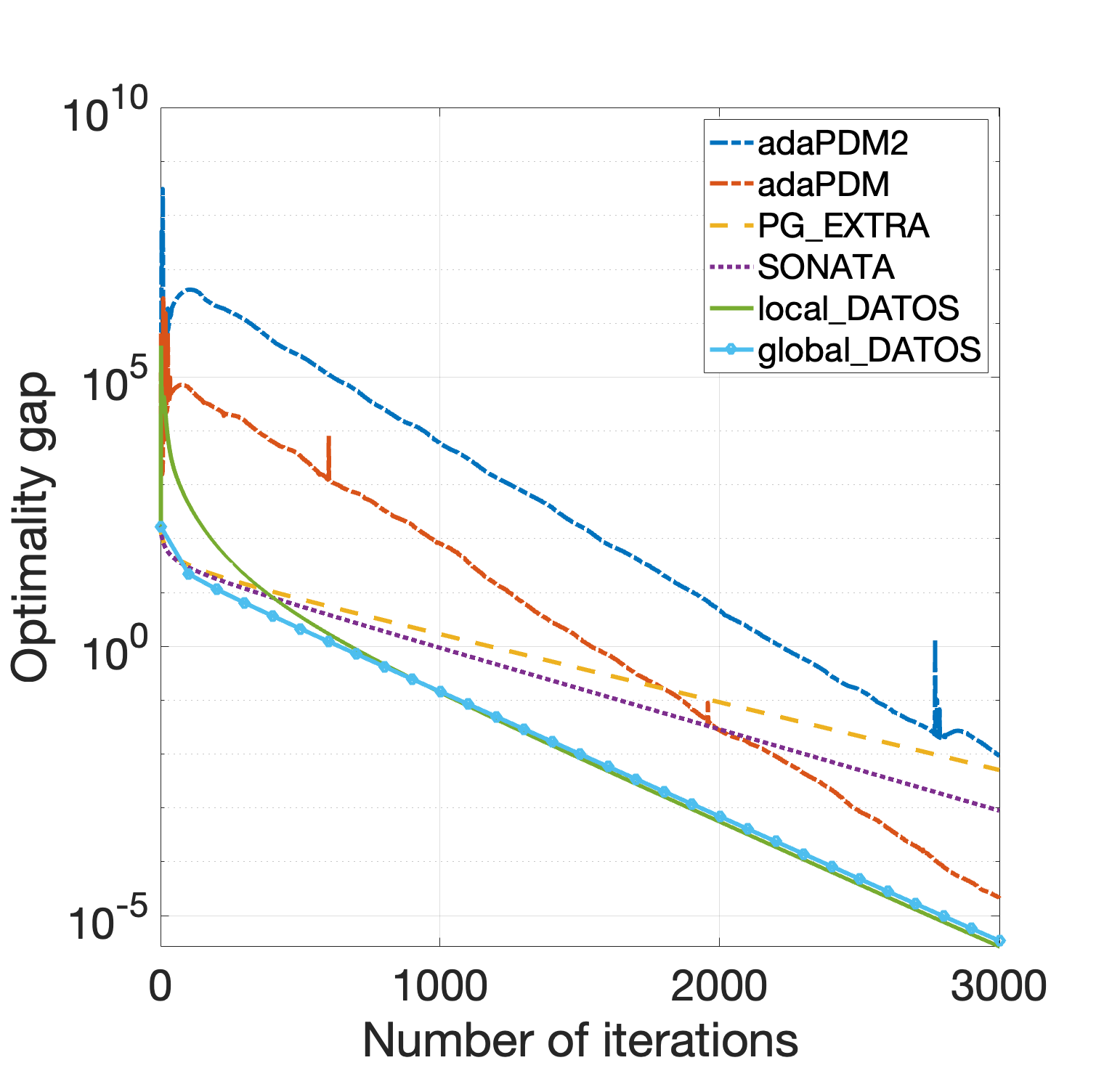}
\caption{Linear regression with elastic net regularization: ${\|\X^k-\X^*\|^2}$  v.s. \# iterations. Comparison of PG-EXTRA, SONATA, adaPDM, adaPDM2, global\_DATOS and local\_DATOS on Erdos-Renyi graphs with different edge-probability:  $p=0.1$  (left);   $p=0.5$ (middle); and  $p=0.9$ (right).}\vspace{-0.45cm} 
\label{fig:general scvx}
\end{figure*}
\begin{theorem}
    \label{thm:convergence local_DATOS}
Consider Algorithm~\ref{alg:DATOS_local} under    Assumptions~\ref{ass:function}-\ref{ass:W}. Further assume that   sequence $\{\T_{\A}^k,\Sub^k\}_k$ is bounded. 
Then, 
   for any given saddle point $(\T^*,\Sub^*)$ of $\mathcal L$, there exists a bounded set   $\mathcal{B}_{T}\times\mathcal{B}_{S}$ such that (i) $(\T^*,\Sub^*)\in \mathcal{B}_{T}\times\mathcal{B}_{S}$, and $(\T_{\A}^k,\Sub^k)\in \texttt{int}(\mathcal{B}_{T}\times\mathcal{B}_{S})$, for all $k\geq 0$; and (ii) for some finite $K\geq 0$, it holds   
  \begin{equation}
        \label{eq:sublinea DATOS_II}
\mathcal{G}_{\mathcal{B}_T\times \mathcal{B}_S}(\overline{\T}_{\A}^k,\overline{\Sub}^k)\leq\frac{1}{k}\cdot\frac{C_K+8R_K^2}{ \min\left(\alpha^{-1},\delta/(2{L}^k)\right)},\quad \forall k\geq K,\end{equation}
where $C_K:=2\sum_{t=1}^K\alpha^{t-1}\mathcal{G}_{\mathcal{B}_T\times\mathcal{B}_S}(\T_{\A}^t,\Sub^t)$ and $$R_K^2:=\|\T_{\B}^{K}-\T^*\|^2+(\alpha^{K})^2\|\Sub^{K}-\Sub^*\|^2.$$

\end{theorem}

      Theorem~\ref{thm:convergence local_DATOS} certifies convergence of   Algorithm~\ref{alg:DATOS} at sublinar rate.  The primary distinction between the local-min and global-min consensus procedures is that the local-min does not guarantee a monotonically decreasing merit function during the initial $K$ iterations. However, as demonstrated by the numerical results in the following section, in practice, Algorithm 3 performs comparably to Algorithm 1. 
 
Convergence of Algorithm 3 relies on the boundedness of the generated iterates. This  is trivially ensured, for instance, when the loss functions are globally smooth or when 
 $r$ is the indicator function of a compact (convex)  set. Additional scenarios and algorithmic variants that guarantee boundedness are  discussed in the extended version of this paper  \cite{chen-et_al-SIOPT25}.
 
\section{Numerical Results}\label{sec:simulations}
This section presents some preliminary numerical results, comparing Algorithms~\ref{alg:DATOS} and~\ref{alg:DATOS_local} against several benchmarks. Specifically, we consider  the {\it nonadaptive} decentralized algorithms SONATA~\cite{sun2019distributed} and PG-EXTRA~\cite{shi2015proximal} as well as decentralized adaptations of the adaptive centralized method adaPDM~\cite{Latafat_23b}.   SONATA and EXTRA require full knowledge of network and optimization parameters; for these methods, we perform a manual grid-search to identify the optimal stepsize ensuring fast convergence. Regarding adaPDM, we emphasize that it is not entirely   parameter-free, since it requires the knowledge of the {\it global} network-related quantity  $\|I-W\|$,  unavailable in practice.  For comparison, we simulate two versions of adaPDM: \textbf{(i)}  one assuming exact knowledge of $\|I-W\|\leq 2$, termed  \texttt{adaPDM}, and \textbf{(ii)} another using the conservative upper bound $2$, resulting in a   network-agnostic variant termed  \texttt{adaPDM2}.  In both cases, we manually tune the parameter $t$ in  adaPDM and adaPDM2 (as in \cite{Latafat_23b}) to obtain the best possible practical convergence behavior.  We anticipate that adaPDM and adaPDM2 are very sensitive to the choice of $t$.

 We simulate  Erdos-Renyi graphs, with $m=20$ agents and  edge-probability of $p=0.1$, $p=0.5$, and $p=0.9$. The gossip weights   used in all the algorithms are the Metropolis-Hasting   (see, e.g., \cite{Nedic_Olshevsky_Rabbat2018}). For the implementation of Algorithm~\ref{alg:DATOS} and~\ref{alg:DATOS_local}, we initialize $\X^0$ and $\Sub^0$ randomly and set $\alpha^{-1}=10$, $\delta=0.9$, and $c=1/3$. 
\label{sec:numerical}
\subsection{Logistic regression with $\ell_1$-regularization}
Consider  the decentralized logistic regression problem with $\ell_1$-regularization, which is an instance of \eqref{eq:problem}, with 
\begin{equation*}\label{eq:logistic-scvx}f_i(x)=\frac{1}{n}\sum_{j=1}^{n}\log(1+\exp(-b_{ij}\cdot\langle x,a_{ij}\rangle)),\quad r(x)=\lambda\|x\|_1,\end{equation*}
where  $a_{ij}\in \mathbb{R}^{d}$ and $b_{ij}\in\{-1,1\}$. The data set  $\{(a_{ij},b_{ij})\}_{j=1}^n$ is  owned   by agent $i$. We use the MNIST dataset from LIBSVM~\cite{chih2libsvm}, taking  the first $N=6000$   samples   (hence $n=300$). The feature dimension is $d=784$. We set the $\ell_1$ regularization parameter $\lambda=10^{-5}$. 

Figure~\ref{fig:general cvx}   plots the optimality gap  $({1}/{m})\sum_{i=1}^mu(x_i)-u^*$ versus the number of iterations, achieved by all the algorihtms,  where $u$ is the objective function in~\eqref{eq:problem} and $u^*$ is its minimal value. 
The figures clearly show  that both  proposed methods consistently outperforms SONATA and PG-EXTRA, with the advantage that do not require any user's intervention for the tuning of the stepsize. Moreover, our methods significantly  outperforms also adaPDM in both of its implementations. Furthermore, the performance difference between Algorithm 3, which incorporates local-min consensus, and Algorithm 1--using global-min consensus,  is remarkably negligible especially in    networks with good connectivity.

\subsection{Maximum likelihood estimate of the covariance matrix}
Consider the decentralized estimation of the inverse of a covariance matrix, which is an instance of \eqref{eq:problem}, with 
\begin{equation*}\label{eq:logistic-scvx}f_i(X)=-n(\log(\texttt{det}(X)))-\texttt{trace}(XY_i),\end{equation*}
and $r(X)=\delta_C(X)$, where $C=\{X\in\mathbb{S}^{d}_{++}: a I\preceq X\preceq bI\}$, for some $0<a\leq b$, and $Y_i=\frac{1}{n}\sum_{j=1}^n y_j^i (y_j^i)^{\top}$ for $\{y_j^i\}$ being the set of local samples at agent $i$'s site, drawn by a Gaussian distribution with covariance matrix $\Sigma\in\mathbb{S}^d_{++}$. Here we take $n=100$ and $d=5$. Notice that this is a problem where the loss function $f$ is only {\it locally} smooth. Hence, SONATA and PG-EXTRA do not have theoretical convergence guarantees. We manually fine-tune their stepsizes for the best stable convergence behavior.

Figure~\ref{fig:general local smooth}   plots the optimality gap  $({1}/{m})\sum_{i=1}^mu(x_i)-u^*$ versus the number of iterations, achieved by all the algorihtms,  where $u$ is the objective function in~\eqref{eq:problem} and $u^*$ is its minimal value calculated by centralized proximal gradient method with linesearch, within   the tolerance of $10^{-30}$. 
The figures clearly illustrate that both proposed methods consistently outperform SONATA and PG-EXTRA, which require conservative stepsize selections to ensure stable global convergence. Additionally, our methods surpass adaPDM in both considered variants. Notice that the iterates generated by Algorithm 3 remain bounded, due to the   compactness enforced  by  $r(X)$.

\subsection{Linear regression with elastic net regularization}
We also report experiments solving  a strongly convex non-smooth instance of Problem~\eqref{eq:problem},  with 
\begin{equation}
    f_i(x)=\frac{1}{n}\|A_ix-b_i\|^2+\frac{\gamma_i}{2}\|x\|^2,\quad r(x)=\lambda\|x\|_1,
\end{equation}
 where $(A_i,b_i)\in\mathbb{R}^{n\times d}\times \mathbb{R}^n$ are the   data  owned  by agent $i$. The elements of $A_i$, $b_i$ are independently sampled from the
standard normal distribution. Here,   $n=20$ (hence $N=400$) and $d=500$. We set $\gamma_i=0.1+(i-1)\times 0.1$ and $\lambda=10^{-5}$, so that the smoothness parameter $L_i$ are different among each local function and the condition number of   $f$ is $\kappa\approx 82.62$.

 Figure~\ref{fig:general scvx} plots the optimality gap measured by $\|\X^k-\X^*\|^2$ versus the number of iterations. Similar to the convex scenario, both proposed algorithms clearly outperform existing decentralized benchmarks, particularly in poorly connected networks. Moreover, the results suggest that our methods exhibit (possibly asymptotic) linear convergence rates--an observation not captured by the  theoretical analysis presented here. Theoretical support for this empirical behavior is provided in the extended version of the paper \cite{chen-et_al-SIOPT25}.

\bibliographystyle{ieeetr-diss}
\bibliography{Bib}

\end{document}